\documentclass[conference]{IEEEtran}
\IEEEoverridecommandlockouts
\usepackage{cite}
\usepackage{amsmath,amssymb,amsfonts}
\usepackage{algorithmic}
\usepackage{graphicx}
\usepackage{textcomp}
\usepackage{xcolor}
\def\BibTeX{{\rm B\kern-.05em{\sc i\kern-.025em b}\kern-.08em
    T\kern-.1667em\lower.7ex\hbox{E}\kern-.125emнX}}

\usepackage{amsmath}
\usepackage{amssymb}
\usepackage{amsfonts}
\usepackage{enumerate}
\usepackage{euscript}
\usepackage{empheq}
\usepackage{mathtools}
\usepackage{tikz}
\usepackage{pgfplots}
\usepackage{euscript}
\usepackage{balance}
\usepackage{cite}
\usepackage{xfrac}

\usetikzlibrary{arrows}
\usetikzlibrary{shapes}
\usetikzlibrary{patterns}
\usepgfplotslibrary{fillbetween}
\usetikzlibrary{arrows.meta}

\newcommand{\minimize}[1]{\underset{{#1}}{\text{min}}}
\newcommand{\maximize}[1]{\underset{{#1}}{\text{max}}}
\newcommand{\fat}[1]{\boldsymbol{#1}}

\newcommand{\braceit}[1]{\left({#1}\right)}
\newcommand{\st}{\text{s.t.}}
\newcommand{\set}{\EuScript}
\usepackage{mathtools}
\DeclarePairedDelimiterX{\norm}[1]{\lVert}{\rVert}{#1}

\newcommand\mydots{\hbox to 0.75em{.\hss.\hss.}}

\tikzset{myarr/.style={{Circle[black,length=4pt]}-{Circle[black,length=4pt]},shorten <=-2.5pt,shorten >=-2.5pt}}

\DeclareMathSymbol{\shortminus}{\mathbin}{AMSa}{"39}
\newcommand{\mysum}[1]{\underset{{{#1}}}{\textstyle\sum}}

\definecolor{red}{HTML}{e05a87} 
\definecolor{blue}{HTML}{035096}

\usepackage[bookmarks=false]{hyperref}

\begin{document}
\title{
Chance-Constrained Equilibrium in Electricity Markets With Asymmetric Forecasts
}

\author{\IEEEauthorblockN{Vladimir Dvorkin, Jalal Kazempour, Pierre Pinson}
\IEEEauthorblockA{\textit{Department of Electrical Engineering} \\
\textit{Technical University of Denmark}\\
Lyngby, Denmark \\
\{vladvo, seykaz, ppin\}@elektro.dtu.dk}
}

\maketitle

\begin{abstract}
We develop a stochastic equilibrium model for an electricity market with asymmetric renewable energy forecasts. In our setting, market participants optimize their profits using public information about a conditional expectation of energy production but use private information about the forecast error distribution. This information is given in the form of samples and incorporated into profit-maximizing optimizations of market participants through chance constraints. We model information asymmetry by varying the sample size of participants' private information. We show that with more information available, the equilibrium gradually converges to the ideal solution provided by the perfect information scenario.  Under information scarcity, however, we show that the market converges to the ideal equilibrium if participants are to infer the forecast error distribution  from  the  statistical  properties of the data at hand or share their private forecasts. 
\end{abstract}



\section{Introduction}

The majority of electricity markets with high penetrations of renewable energy sources clear several trading floors, e.g., day-ahead and in real-time, to offset  potential imbalances induced by renewable forecast errors. Given a forecast of renewable generation, the day-ahead stage determines the optimal allocation of energy and reserves to offset any forecast deviation at the real-time stage. To enable reliable and cost-efficient operations, one of the existing suggestions in the technical literature is to optimize the day-ahead decisions using stochastic programming \cite{morales2009economic}. By taking a probabilistic forecast as input, either in the form of discrete scenarios \cite{juanmi} or moments of forecast error distribution \cite{dvorkin2019chance}, stochastic models produce robust day-ahead decisions. 

A common assumption in the literature is that all market participants use identical information about the uncertainty distribution of renewable generation. We refer to this situation as forecast symmetry. However, market participants may use various forecasting tools of different quality or source their forecasts from different providers. Market participants often treat their forecasts as private data and have no means of evaluating the benefits from sharing their forecasts. Moreover, even with identical forecasts, market participants with heterogeneous risk attitudes \cite{bertsimas2009constructing,gerard2018risk} or irrational preferences \cite{tversky1992advances} utilize the available data differently. Hence, we relax the assumption of forecast symmetry and explore the impacts of \textit{asymmetric} renewable forecasts on electricity market outcomes.

There are a few works addressing forecast asymmetry and its impacts on electricity market outcomes.  Using a scenario-based stochastic programming and a game-theoretic analysis, \cite{dvorkin2019electricity} illustrates how the social welfare in competitive electricity markets varies as a function of the level of forecast asymmetry among market participants. In oligopolistic setting, \cite{lazaros} analyzes the impact of public dissemination of aggregate renewable power forecast on market outcomes.

In this work, we study the impacts of renewable forecast asymmetry among market participants on the social and individual market outcomes. To assess these impacts, we build a stochastic equilibrium model, including a set of profit-maximizing optimization problems, one per market participant, coupled by power balance conditions. In our model, each market participant optimizes its expected profit, while using its own private information about forecast error distribution. This information is obtained from independent providers and given in a form of \textit{samples}. To model the forecast asymmetry, we vary the sample size of private forecast datasets. To incorporate private forecasts into the stochastic equilibrium model, we develop a \textit{chance-constrained} optimization for profit-maximization problem of each market participant. Using the sample representation of forecast errors, the objective function and the feasible region of market participant problem are made conditional on the private forecast. 

Using a stylized case study, we conduct three experiments. First, we explore the market implications of forecast asymmetry and illustrate that system reliability and operating cost may significantly improve even with a marginal increase in the sample size of private forecast datasets. We also show that the system reliability converges to a desirable level comparatively faster than the operating cost. Then, we study two approaches to enhance the market operation under data scarcity. We introduce a case where market participants are able to infer the stochastic process distribution from the statistical properties learned from the data at hand. Through learning, the reliability significantly improves even if the size of dataset provided by forecast providers is relatively small. Third, we show that there might be some circumstances under which market participants  have strong incentives for sharing their private forecasts as it improves not only the overall market performance, but also their individual profit outcomes.

The rest of this paper is outlined as follows. Section II details the transition from a centralized stochastic model with symmetric forecasts to a stochastic equilibrium model with forecast asymmetry. Section III incorporates private forecasts into individual optimization problems. Section IV streamlines the decentralized algorithm to compute equilibrium solution. Section V provides the results of the three experiments on a stylized system. Finally, Section VI concludes.

\section{Market clearing: From a centralized optimization to an equilibrium model}
\subsection{Preliminaries}
We consider a day-ahead electricity market with a high share of renewable energy production. This market is cleared 12--36 hours ahead of real time to meet total load $L$ and offset any imbalance induced by uncertain renewable energy production in a reliable and cost-efficient manner. The expected renewable energy production is given by a point forecast $w^{f}$ assumed to be a public information. The real-time deviation from $w^{f}$ is given by a zero-mean random forecast error denoted by $\fat{w}$. The distribution of $\fat{w}$ is considered to be a private information for each market participant, and thereby, market participants may have asymmetric information about the random forecast error $\fat{w}$. Note that we make no assumption on the type of distribution of $\fat{w}$. The set $\set{G}$ solely includes conventional power producers. Following an affine policy for reserve allocation \cite{joe}, the eventual production of each conventional producer $i\in\set{G}$ in real-time, i.e., $\fat{p}_{i}$, is considered as an affine function of forecast error, i.e., 
\begin{align}
    \fat{p}_{i}=p_{i} -\alpha_{i}\fat{w}, \; \alpha_{i} \geqslant 0, \quad i\in\set{G}, \label{control_policy}
\end{align}
where $p_{i}$ is the nominal dispatch in the day-ahead stage and $\alpha_{i}$ is a portion of renewable power deviation adjusted by producer $i$, the so-called \textit{participation factor}. Both $p_{i}$ and $\alpha_{i}$ are day-ahead stage decision variables. By using \eqref{control_policy},  producers co-optimize their day-ahead dispatch decision and participation in real-time adjustment with respect to the realization of the forecast error.  To ensure that the entire deviation from the forecast is accommodated, we enforce the adjustment balance as  $\sum_{i\in\set{G}}\alpha_{i}=1$. Each producer $i$ outputs within its operational limits $[\underline{p}_{i},\overline{p}_{i}\big]$ and adjusts its production in real-time by at most $\overline{r}_{i}$. The production cost of each producer is quadratic with the first- and second-order coefficients $c_{1i}$ and $c_{2i}$, respectively. 
\subsection{Centralized optimization model}
Under an assumption of identical forecast for all market participants, the market is cleared in a centralized manner using a  stochastic optimization model as
\begin{subequations}
\begingroup
\allowdisplaybreaks
\label{CC_dispatch}
\begin{align}
    \minimize{p,\alpha\geqslant0}\quad&\mathbb{E}\bigg[\mysum{i\in\set{G}}\big(c_{2i}\fat{p}_{i}^2 + c_{1i}\fat{p}_{i}\big)\bigg]\label{CC_dispatch_obj}\\
    \st \quad\mathbb{P}&\begin{pmatrix*}[l]&\hspace{-1em}p_{i} - \alpha_{i}\fat{w} \leqslant \overline{p}_{i},\;\forall i\in\set{G}\\
    &\hspace{-1em}p_{i} - \alpha_{i}\fat{w} \geqslant \underline{p}_{i},\;\forall i\in\set{G}\\
    &\hspace{-1em}\alpha_{i}\fat{w} \leqslant \overline{r}_{i},\;\forall i\in\set{G} \\
    &\hspace{-1em}\alpha_{i}\fat{w} \geqslant -\overline{r}_{i},\;\forall i\in\set{G}\end{pmatrix*}  \geqslant 1-\varepsilon,\label{prob_con}\\
    &\mysum{i\in\set{G}}p_{i} + w^{f} - L = 0 : \lambda^{\text{e}},\label{energy_bal}\\
    &\mysum{i\in\set{G}} \alpha_{i} = 1 : \lambda^{\text{r}},\label{res_bal}
\end{align}
\endgroup
\end{subequations}
where $\mathbb{E}(\cdot)$ is the expectation operator and $\mathbb{P}(\cdot)$ is the probability operator, whereas the random forecast error $\fat{w}$ follows a distribution with known parameters for all market participants. The objective function \eqref{CC_dispatch_obj} minimizes the expected total operating cost of the system. The  chance constraint \eqref{prob_con} is enforced to ensure that constraints of all producers hold jointly with a probability at least  $(1-\varepsilon)$.
The maximum allowable constraint violation probability $\varepsilon$ is a measure of the acceptable level of risk exposure and is kept small to ensure a reliable real-time operation. The equality constraint \eqref{energy_bal} enforces the power production and consumption balance in the day-ahead stage. Finally, \eqref{res_bal} ensures allocating sufficient reserve, such that the entire renewable power production imbalance in real-time will be offset. The dual variables of the equality constraints define the energy price $\lambda^{\text{e}}$ and the reserve price $\lambda^{\text{r}}$. 

\subsection{Equilibrium model}
By design, the centralized problem \eqref{CC_dispatch} does not incorporate the private information on forecast error distribution. To model market outcomes with private forecasts, we introduce the following stochastic equilibrium problem: 
\begin{subequations}
\begingroup
\allowdisplaybreaks
\label{eq}
\begin{align}
&\left\{\!
\begin{aligned}
    \maximize{p_{i},\alpha_{i}\geqslant0} \; & \lambda^{\text{e}}p_{i} - \lambda^{\text{r}}\alpha_{i} -\mathbb{E}_{i}\big[c_{2i}\fat{p}_{i}^2 + c_{1i}\fat{p}_{i}\big] \\
    \st \;\mathbb{P}_{i}&\begin{pmatrix*}[l]&\hspace{-1em}p_{i} - \alpha_{i}\fat{w} \leqslant \overline{p}_{i}\\
    &\hspace{-1em}p_{i} - \alpha_{i}\fat{w} \geqslant \underline{p}_{i}\\
    &\hspace{-1em}\alpha_{i}\fat{w} \leqslant \overline{r}_{i} \\
    &\hspace{-1em} \alpha_{i}\fat{w} \geqslant -\overline{r}_{i}\end{pmatrix*}  \geqslant 1-\varepsilon_{i}
\end{aligned}
\right\},\;\forall i \in\set{G},\label{eq_producer_problem}\\
&\maximize{\lambda^{\text{e}}}\quad\lambda^{\text{e}}\big[\mysum{i\in\set{G}}p_{i} + w^{f} - L\big],\label{eq_price_setter_energy}\\
&\maximize{\lambda^{\text{r}}}\quad\lambda^{\text{r}}\big[\mysum{i\in\set{G}}\alpha_{i} - 1\big],\label{eq_price_setter_regulation}
\end{align}
\endgroup
\end{subequations}
where each conventional producer $i\in\set{G}$ maximizes its expected profit in \eqref{eq_producer_problem} for given energy and reserve prices. Besides, for given values of $p_{i}$ and $\alpha_{i}$, the unconstrained problems \eqref{eq_price_setter_energy} and \eqref{eq_price_setter_regulation} set the energy and reserve prices. One can observe that the Karush--Kuhn--Tucker conditions of \eqref{eq_price_setter_energy} and \eqref{eq_price_setter_regulation} respectively yield the balancing constraints \eqref{energy_bal} and \eqref{res_bal}. Unlike the centralized optimization \eqref{CC_dispatch}, $\mathbb{E}_{i}(\cdot),$ $\mathbb{P}_{i}(\cdot),$ and $\varepsilon_{i}$ in \eqref{eq_producer_problem} are indexed by $i$, indicating that each conventional producer incorporates its own information about forecast error distribution. With asymmetric information, the solution of \eqref{eq} does not amount to that of \eqref{CC_dispatch}, and the corresponding market implications are the main focus of this paper. 



\section{Incorporation of private forecasts and chance constraints reformulation}
We model the private information of each conventional producer $i\in\set{G}$ on forecast error distribution $\fat{w}$ by a finite set of samples, i.e., $\mathcal{D}_{i}:=\{w_{i1},\hspace{-0.2em}...,w_{is}\},\;\forall i\in\set{G}. $
This sample representation enables a distribution-free reformulation of individual optimization problems \eqref{eq_producer_problem}. 

We first start with reformulating the objective function of \eqref{eq_producer_problem}. The expected production cost for each conventional producer expresses as 
\begin{align}
    \mathbb{E}_{i}\big[c_{2i}\fat{p}_{i}^2 + c_{1i}\fat{p}_{i}\big]&=\mathbb{E}_{i}\big[c_{2i}(p_{i} -\alpha_{i}\fat{w})^2 + c_{1i}(p_{i} -\alpha_{i}\fat{w})\big]\nonumber\\
    &= c_{2i}p_{i}^2 + c_{1i}p_{i} + \mathbb{E}_{i}\big[c_{2i}(\alpha_{i}\fat{w})^2\big]\nonumber\\
    &= c_{2i}p_{i}^2 + c_{1i}p_{i} + c_{2i}(\alpha_{i}\hat{\sigma}_{i})^2,\label{cost_reformulation}
\end{align}
where $\hat{\sigma}_{i}^2$ denotes the variance of the forecast error distribution obtained from the sample-based dataset $\mathcal{D}_{i}.$ 

With the sample representation of random variable $\fat{w}$, a straightforward methodology to reformulate the joint chance constraint in \eqref{eq_producer_problem} is to enforce its entries on each sample of the forecast error distribution \cite{alamo2010sample}. However, this methodology is computationally expensive as the size of dataset required for achieving a probabilistic performance guarantee grows in the number of decision variables \cite[Theorem 4]{alamo2010sample}. To ease the computational burden, we refer to the work in \cite{margellos2014road} that offers a compromised solution between sample approximation and robust optimization. In \cite{margellos2014road}, the constraints are enforced over the vertices (bounds for scalar uncertainties) of the minimum volume of the hyper-rectangular uncertainty set. This allows to re-define the requirement for the size of dataset as in \cite[Equations (7)]{margellos2014road}, which is independent from the number of decision variables. Instead, the  size of sample-based dataset required for the chance constraint to hold grows in the number of uncertainty sources, which amounts to one in our case. 

In line with \cite{margellos2014road}, we consider that conventional producers bound the support of forecast error distribution as $\underline{w}_{i}:=\min\{w_{i1},\hspace{-0.2em}...,w_{is}\}$ and $\overline{w}_{i}:=\max\{w_{i1},\hspace{-0.2em}...,w_{is}\}$, yielding $\underline{w}_{i}\leqslant0\leqslant\overline{w}_{i}$. Now, for each conventional producer, the joint chance constraint \eqref{eq_producer_problem} is approximated through a set of inequalities
\begin{subequations}
\label{reformulated_cc}
\begin{align}
    & p_{i} - \alpha_{i} \underline{w}_{i}\leqslant \overline{p}_{i},\quad \alpha_{i} \underline{w}_{i} \geqslant -\overline{r}_{i},\\ 
    & p_{i} - \alpha_{i} \overline{w}_{i}\geqslant \underline{p}_{i},\quad \alpha_{i} \overline{w}_{i} \leqslant \overline{r}_{i}.
\end{align}
\end{subequations}
These  inequalities intuitively exhibit the direct impact of sample-based dataset $\mathcal{D}_{i}$ on reserve margins required to provide adjustment in real time. With smaller sample support $[\underline{w}_{i},\overline{w}_{i}]$, the producers reserve less capacity to offset the forecast error in real time. 

\begin{figure}
\centering
\resizebox{0.48\textwidth}{!}{%
\begin{tikzpicture}
\def\normaltwo{\x,{1.25*1/exp(((\x-5)^2)/2)}}
\def\y{7} 
\def\z{3} 
\fill [fill=gray!35] (\z,0) -- plot[domain=\z:\y,samples=100] (\normaltwo) -- ({\y},0) -- cycle;
\draw[color=black,line width = 0.15mm, domain=\z:\y,samples=100] plot (\normaltwo) node[right] {};
\draw[black!50] ({7},{0.17}) -- ({7},0) node[above,font=\scriptsize,yshift=5] {\textcolor{black}{$\overline{w}$}};
\draw[black!50] ({3},{0.17}) -- ({3},0) node[above,font=\scriptsize,yshift=5] {\textcolor{black}{$\underline{w}$}};
\draw[black!50,dashed] (1.3,0) -- (3,0);
\draw[black!50,dashed] (7,0) -- (8.7,0);
\draw[->,>=stealth,black] (2.25,0) -- (7.75,0);
\draw[black] (7.75,0) node[above,font=\scriptsize,xshift=0,yshift=0] {$\fat{w}$};
\draw[->,>=stealth,black] (5,0) -- node[rotate=90,font=\scriptsize,xshift=17,yshift=4.5] {pdf} (5,1.75);
\draw[black!50,dashed] (3,-0.) -- (3,-2.7) node[right] {};
\draw[black!50,dashed] (7,-0.) -- (7,-2.7) node[right] {};
\draw[black!50,dashed] (1.3,-0.9) -- (8.7,-0.9) node[right] {};
\draw[black!50,dashed] (3.75,-0.35) -- (3.75,0) node[right] {};
\draw[black!50,dashed] (6,-0.55) -- (6,0) node[right] {};
\draw[blue] (3.75,-0.35) -- (5.5,-0.35) node[right] {};
\node[mark size=1pt,color=blue] at (3.75,-0.35) {\pgfuseplotmark{*}};
\node[mark size=1pt,color=blue] at (5.5,-0.35) {\pgfuseplotmark{*}};
\draw[black] (3.75,-0.35) node[above,font=\scriptsize,yshift=-0.75,xshift=6] {$\underline{w}_{1}$};
\draw[black] (5.5,-0.35) node[above,font=\scriptsize,yshift=-0.75] {$\overline{w}_{1}$};
\draw[red] (4.5,-0.55) -- (6,-0.55) node[right] {};
\node[mark size=1pt,color=red] at (4.5,-0.55) {\pgfuseplotmark{*}};
\node[mark size=1pt,color=red] at (6,-0.55) {\pgfuseplotmark{*}};
\draw[black] (4.5,-0.55) node[below,font=\scriptsize,yshift=0.75] {$\underline{w}_{i}$};
\draw[black] (6,-0.55) node[below,font=\scriptsize,yshift=0.75] {$\overline{w}_{i}$};
\draw[draw] (2.1,-0.08) node[below,align=center,font=\scriptsize] {small \\ reliability};
\draw[draw] (7.9,-0.08) node[below,align=center,font=\scriptsize] {high cost \\ of operation};
\node[draw=none,rotate=-15,black!70] at (5,-0.45) {...};
\draw[black!50,dashed] (1.3,-1.8) -- (8.7,-1.8) node[right] {};
\draw[blue] (3.0,-1.25) -- (5.5,-1.25) node[right] {};
\node[mark size=1pt,color=blue] at (3.0,-1.25) {\pgfuseplotmark{*}};
\node[mark size=1pt,color=blue] at (5.5,-1.25) {\pgfuseplotmark{*}};
\draw[black] (3.0,-1.25) node[above,font=\scriptsize,xshift=6,yshift=-0.75] {$\underline{w}_{1}$};
\draw[black] (5.5,-1.25) node[above,font=\scriptsize,xshift=0,yshift=-0.75] {$\overline{w}_{1}$};
\draw[red] (4.5,-1.45) -- (7,-1.45) node[right] {};
\node[mark size=1pt,color=red] at (4.5,-1.45) {\pgfuseplotmark{*}};
\node[mark size=1pt,color=red] at (7,-1.45) {\pgfuseplotmark{*}};
\draw[black] (4.5,-1.45) node[below,font=\scriptsize,yshift=0.75] {$\underline{w}_{i}$};
\draw[black] (7,-1.45) node[below,font=\scriptsize,xshift=-6,yshift=0.75] {$\overline{w}_{i}$};
\draw[draw] (2.1,-0.95) node[below,align=center,font=\scriptsize] {high \\ reliability};
\draw[draw] (7.9,-0.95) node[below,align=center,font=\scriptsize] {medium cost \\ of operation};
\node[draw=none,rotate=-15,black!70] at (5,-1.35) {...};
\draw[black!50,dashed] (1.3,-2.7) -- (8.7,-2.7) node[right] {};
\draw[blue] (3.0,-2.15) -- (7,-2.15) node[right] {};
\node[mark size=1pt,color=blue] at (3.0,-2.15) {\pgfuseplotmark{*}};
\node[mark size=1pt,color=blue] at (7,-2.15) {\pgfuseplotmark{*}};
\draw[black] (3.0,-2.15) node[above,font=\scriptsize,xshift=6,yshift=-0.75] {$\underline{w}_{1}$};
\draw[black] (7,-2.15) node[above,font=\scriptsize,xshift=-6,yshift=-0.75] {$\overline{w}_{1}$};
\draw[red] (3,-2.35) -- (7,-2.35) node[right] {};
\node[mark size=1pt,color=red] at (3,-2.35) {\pgfuseplotmark{*}};
\node[mark size=1pt,color=red] at (7,-2.35) {\pgfuseplotmark{*}};
\draw[black] (3,-2.35) node[below,font=\scriptsize,xshift=6,yshift=0.75] {$\underline{w}_{i}$};
\draw[black] (7,-2.35) node[below,font=\scriptsize,xshift=-6,yshift=0.75] {$\overline{w}_{i}$};
\draw[draw] (2.1,-1.85) node[below,align=center,font=\scriptsize] {high \\ reliability};
\draw[draw] (7.9,-1.85) node[below,align=center,font=\scriptsize] {small cost \\ of operation};
\node[draw=none,rotate=-15,black!70] at (5,-2.25) {...};
\end{tikzpicture}
}
\caption{Relation between the size of sample-based dataset of conventional producers and the market operation.}
\label{example}
\end{figure}
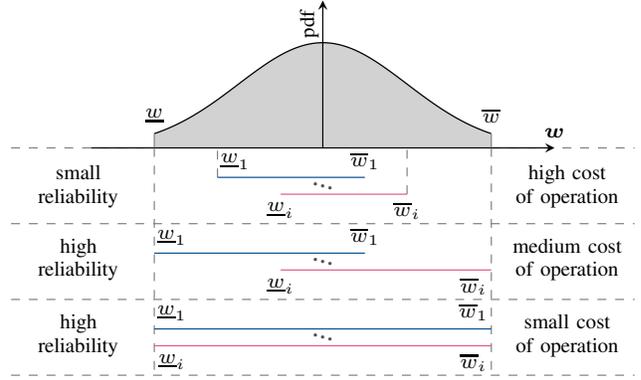

Notice, in what follows, we do not specify $\varepsilon_{i}$ for each producer $i$ required for sample generation in \cite{margellos2014road}. Instead, we directly analyze the dependency of the system-wide reliability on the size of sample-based dataset of conventional producers. We explain this rationale with an illustration in Fig. \ref{example}. With a relatively small size of sample-based dataset, there is a risk of not capturing the entire range of forecast error realizations, eventually leading to a small system reliability. As a result, expensive extreme balancing actions, such as renewable power spillage or load curtailment, would be required to restore the real-time balance. With an increasing size of sample-based dataset, the union of datasets of producers may suffice to keep the system in balance with high probability, which improves the system reliability. 
With a sufficiently large size of sample-based dataset, all producers are available to balance any realization of the forecast error. In this scenario, the dispatch is solely driven by the marginal cost of generation, thus improving the overall operating cost.

\section{Equilibrium computation}

\begin{figure}
\tikzstyle{box} = [rectangle, rounded corners = 10, minimum width=100, minimum height=30,text centered, draw=black, fill=white!1,line width=0.3mm]
\tikzstyle{infobox} = [rectangle, rounded corners = 10, minimum width=50, minimum height=20,text centered, draw=black, fill=white!1,line width=0.3mm]
\centering
\resizebox{0.44\textwidth}{!}{%
\begin{tikzpicture}[node distance=100]
\node [align=center] (Generator1) [box] {Update of producer 1\\[0.5mm]  Update $p_{1},\alpha_{1}$ using \eqref{producerupdate}};
\node [align=center] (GeneratorI) [box,xshift=202.5] {Update of producer $i$\\[0.5mm]  Update $p_{i},\alpha_{i}$ using \eqref{producerupdate}};
\node [align=center] (Energy_price) [box,xshift=5cm,yshift=2.15cm,xshift=-40] {Energy price update\\[0.5mm]  Update $\lambda^{\text{e}}$ using \eqref{energypriceupdate}};
\node [align=center] (Reserve_price) [box,xshift=5cm,yshift=-2.15cm,xshift=-40] {Reserve price update\\[0.5mm]  Update $\lambda^{\text{r}}$ using \eqref{reservepriceupdate}};
\node [align=center] (others) [below of = Energy_price, yshift=1.25cm] {$\boldsymbol{\dots}$};
\path[] 
(Generator1) edge[bend left,->,>=stealth,line width = 0.35mm] node [left, sloped, above] {$p_{1}$}(Energy_price)
(GeneratorI) edge[bend right,->,>=stealth,line width = 0.35mm] node [left, sloped, above] {$p_{i}$}(Energy_price)
(Energy_price) edge[bend left,->,>=stealth,line width = 0.35mm] node [left, sloped, below] {$\lambda^{\text{e}}$}(Generator1)
(Energy_price) edge[bend right,->,>=stealth,line width = 0.35mm] node [left, sloped, below] {$\lambda^{\text{e}}$}(GeneratorI);
\path[] 
(Generator1) edge[bend right,->,>=stealth,line width = 0.35mm] node [left, sloped, below] {$\alpha_{1}$}(Reserve_price)
(GeneratorI) edge[bend left,->,>=stealth,line width = 0.35mm] node [left, sloped, below] {$\alpha_{i}$}(Reserve_price)
(Reserve_price) edge[bend right,->,>=stealth,line width = 0.35mm] node [left, sloped, above] {$\lambda^{\text{r}}$}(Generator1)
(Reserve_price) edge[bend left,->,>=stealth,line width = 0.35mm] node [left, sloped, above] {$\lambda^{\text{r}}$}(GeneratorI);
\end{tikzpicture}
}
\caption{Proposed iterative approach based on Walrasian t\^{a}tonnement to determine the market equilibrium point.}
\label{scheme}
\end{figure}
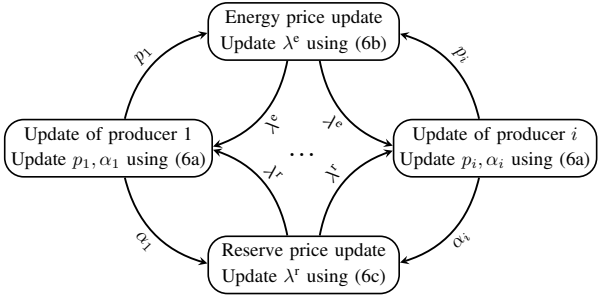

We compute the solution to the equilibrium problem \eqref{eq} using an iterative algorithm schematically depicted in Fig. \ref{scheme}. The algorithm is inspired by Walrasian t\^{a}tonnement \cite{uzawa1960walras}. Each producer updates its nominal dispatch and participation factor based on the corresponding prices for energy and reserve. Upon receiving updates from producers, the price-setting problems adjust the prices. The algorithm formulates as the following iterating procedure: 
\begin{subequations}
\begingroup
\allowdisplaybreaks
\label{algorithm}
\begin{align}
    (p,\alpha)\leftarrow&\;\text{argmin}\;\eqref{eq_producer_problem}\;\text{for given}\;\lambda^{\text{e}},\lambda^{\text{r}},\label{producerupdate}\\
    \lambda^{\text{e}} \leftarrow& \lambda^{\text{e}} - \rho\big[\mysum{i\in\set{G}}p_{i} + w^{f} - L\big],\label{energypriceupdate}\\
    \lambda^{\text{r}} \leftarrow& \lambda^{\text{r}} - \rho\big[\mysum{i\in\set{G}}\alpha_{i}-1\big],\label{reservepriceupdate}
\end{align}
\endgroup
\end{subequations}
where prices in \eqref{energypriceupdate} and \eqref{reservepriceupdate} evolve through iterations along the decent directions $\nabla_{\lambda^{\text{e}}}\eqref{eq_price_setter_energy}$ and $\nabla_{\lambda^{\text{r}}}\eqref{eq_price_setter_regulation}$, respectively, with a suitable step size $\rho\geqslant0$. If the aggregated energy and reserve exceed the corresponding market needs, energy and reserve prices will be reduced to minimize the imbalance. Likewise, these prices increase for any shortage of aggregated energy and reserve supplies. This requires  $\rho$ to be properly selected.

As objective function of each producer in \eqref{eq_producer_problem} is strictly convex in their decision variables, the algorithm in \eqref{algorithm} provably converges to the solution of equilibrium problem \eqref{eq}, provided that a solution exists \cite{falsone2017dual}. 

\section{Numerical experiments}\label{numerical_experiments}
In Section \ref{eval_framework} we introduce an evaluation framework to assess the efficiency of the proposed chance-constrained equilibrium model with forecast asymmetry. In Section \ref{setup} we explain a stylized system setup. In Sections \ref{exp1} to \ref{exp3}, we introduce a set of experiments to (a) analyze the impacts of the number of samples within the dataset of each producer on the overall system performance, (b) analyze the system-wide benefits when producers  learn from their datasets, and (c) show the overall and individual benefits from data sharing.

\subsection{Evaluation framework} \label{eval_framework}
\begin{figure}
\centering
\tikzstyle{mainbox} = [rectangle, rounded corners = 10, minimum width=100, minimum height=30,text centered, draw=black, fill=white!5,line width=0.3mm]
\tikzstyle{producer} = [rectangle, rounded corners = 10, minimum width=30, minimum height=30,text centered, draw=black, fill=white!5,line width=0.3mm]
\tikzstyle{Out_of_sample} = [rectangle, rounded corners = 10, minimum width=100, minimum height=30,text centered, draw=black, fill=white!5,line width=0.3mm]
\resizebox{0.45\textwidth}{!}{%
\begin{tikzpicture}[node distance=62]
\node [align=center] (Algorithm) [mainbox] {Compute optimal set \\
points using \eqref{algorithm}};
\node [align=center] (producer1) [producer, above of = Algorithm] {Data processing \\ by producer 1};
\node [align=center] (Provider1) [producer, left of = producer1, xshift = -45] {Data provider \\ of producer 1};
\node [align=center] (producer2) [producer, below of = Algorithm] {Data processing \\ by producer $i$};
\node [align=center] (Provider2) [producer, left of = producer2, xshift = -45] {Data provider \\ of producer $i$};
\node [align=center] (others) [left of = Algorithm, rotate = 90, yshift = 45] {$\fat{\dots}$};
\node [align=center] (Out_of_sample1) [Out_of_sample, right of = Algorithm, xshift = 65] {Solve re-dispatch\\
problem \eqref{out_of_sample_opt} $\forall s$};
\node [align=center] (Out_of_sample_data) [Out_of_sample, above of = Out_of_sample1] {Out-of-sample \\ data};
\node [align=center] (Results) [producer, below of = Out_of_sample1] {Results};
\draw[->,>=stealth,line width=0.35mm](producer1.south)-- node [sloped, above] {$\hat{\sigma}_{1}^2$}(Algorithm.north);
\draw[->,>=stealth,line width=0.35mm](producer1.south)-- node [sloped, below] {$\underline{w}_{1},\overline{w}_{1}$}(Algorithm.north);
\draw[->,>=stealth,line width=0.35mm](producer2.north)-- node [sloped, above] {$\hat{\sigma}_{i}^2$}(Algorithm.south);
\draw[->,>=stealth,line width=0.35mm](producer2.north)-- node [sloped, below] {$\underline{w}_{i},\overline{w}_{i}$}(Algorithm.south);
\draw[->,>=stealth,line width=0.35mm](Provider1.east)-- node [sloped, above] {$\mathcal{D}_{1}$}(producer1.west);
\draw[->,>=stealth,line width=0.35mm](Provider2.east)-- node [sloped, above] {$\mathcal{D}_{i}$}(producer2.west);
\draw[->,>=stealth,line width=0.35mm](Algorithm.east)-- node [sloped, above] {$p_{i}^{\star}$}(Out_of_sample1.west);
\draw[->,>=stealth,line width=0.35mm](Out_of_sample_data.south)-- node [sloped, above] {$\mathcal{D}$}(Out_of_sample1.north);
\draw[->,>=stealth,line width=0.35mm](Out_of_sample1.south)-- node [sloped, above] {$\nu,\mathcal{C}$}(Results.north);
\end{tikzpicture}
}
\caption{Evaluation framework.}
\label{evaluation_framework}
\end{figure}
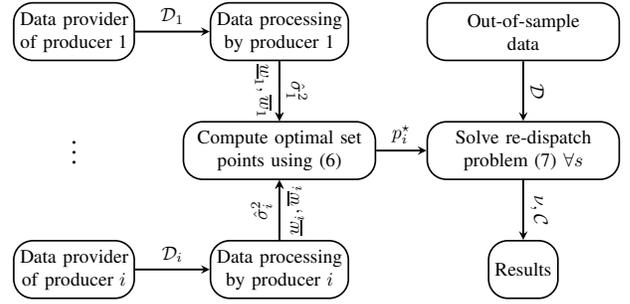

To assess the impacts of forecast asymmetry among conventional producers on the overall market performance, we use the evaluation framework schematically depicted in Fig. \ref{evaluation_framework}. Here, each producer $i$ receives a sample-based dataset $\mathcal{D}_{i}=\{w_{i1},\hspace{-0.2em}...,w_{is}\}$ with $N_{i}^{s}$ samples from its private data provider. In the next step, producers process their own data to estimate the variance $\hat{\sigma}_{i}$ and the support of the forecast error distribution contained between bounds $\underline{w}_{i}$ and $\overline{w}_{i}$. Then, the optimal nominal dispatch $p_{i}^{\star}$ and participation factor $\alpha_{i}^{\star}$ are obtained from algorithm \eqref{algorithm}. To assess the efficiency of day-ahead decisions $p_{i}^{\star}$ in real-time, we make an out-of-sample analysis over the set $\mathcal{D}=\{w_{1},\hspace{-0.2em}...,w_{s}\}$ with $N$ number of samples. We choose $N=3\cdot10^5$. In all experiments, we consider $\mathcal{D}_{i}$ and $\mathcal{D}$ containing independent and identically distributed random samples. 

One of the common methods to assess the constraint violation level is to verify the frequency of feasible solution given constraints \eqref{reformulated_cc} and optimal values $p_{i}^{\star}$. The method assumes that producers are equipped with automatic control and respond to uncertainty realization according to the affine function \eqref{control_policy}. However, this approach overlooks the cost associated with the real-time extreme actions such as renewable power spillage and load curtailment. Adopting the approach from \cite{ordoudis2018energy}, we consider, instead, that the system is capable of re-dispatching power producers closer to real-time operations. Taking into account $w_{s}$ as the renewable power production realized in real-time, each producer adjusts its production in real-time denoted by $r_{i}$ to offset the imbalance in a cost-efficient manner. The corresponding re-dispatch problem in real time formulates as
\begin{subequations}
\begingroup
\allowdisplaybreaks
\label{out_of_sample_opt}
\begin{align}
    \minimize{r_{i},w^{c},\atop l^{s},c^{e}}& \mathcal{C} := \sum_{i\in\set{G}} \braceit{c_{2i}\braceit{p_{i}^{\star} + r_{i}}^2 + c_{1i}\braceit{p_{i}^{\star} + r_{i}}} + c^{e} \label{out_opt_obj}\\
    \st\;\quad&\sum_{i\in\set{G}} r_{i} + l^{s} + w_{s} - w^{c} = 0,\label{out_opt_con1}\\
    &\underline{p}_{i} \leq p_{i}^{\star} + r_{i} \leq \overline{p}_{i}, \quad \forall i\in\set{G},\label{out_opt_con3}\\
    &-\overline{r}_{i} \leq r_{i} \leq \overline{r}_{i}, \quad \forall i\in\set{G},\label{out_opt_con4}\\
    &c^{e}=c^{c}w^{c} + c^{s}l^{s},\label{out_opt_con2}\\
    &0\leq w^{c} \leq w^{f} + w_{s}, \label{out_opt_con5}\\
    &0\leq l^{s} \leq L,\label{out_opt_con6}
\end{align}
\endgroup
\end{subequations}
where the objective function \eqref{out_opt_obj} includes the total generation cost of conventional producers as well as emergency cost $c^{e}$ induced by renewable power spillage $w^{c}$ and load curtailment $l^{s}$, required when the market is in deficit of reserve capacity. The renewable power spillage cost $c^{c}$ and load curtailment cost $c^{s}$ are commonly set high, so that the market prioritizes adjusting the production of flexible conventional producers. Observe, that the generating cost in \eqref{out_opt_obj} is a deterministic variant of \eqref{CC_dispatch_obj} once nominal dispatch for all producers $p_{i}^{\star}$ is obtained. Similarly, constraints of producers in \eqref{out_opt_con3}-\eqref{out_opt_con4} are indeed deterministic variants of \eqref{prob_con}. Whenever either $w^{c}$ or $l^{s}$ is positive, we record the empirical violation $\nu$ of constraint \eqref{prob_con} and corresponding reliability level $(1-\nu)$. Despite original constraint set violation, the problem remains feasible due to additional flexibility from $w^{c}$ or $l^{s}$. Therefore, we can assess operating cost $\mathcal{C}$ for any realization of uncertainty. Notice, that the centralized re-dispatch problem \eqref{out_of_sample_opt} is solely used for the evaluation procedure, and its decentralized counterpart is achieved by the same means of algorithm \eqref{algorithm}. 

As producers receive different samples at each run of the evaluation framework, we make 50 simulation runs for each $N_{i}^{s}$ to increase the statistical significance of the results. 

\subsection{Experimental setup}\label{setup}
We consider a set of two conventional producers $\set{G}\in\{1,2\}$ with maximum capacity $\overline{p}_{i}=\{32,44\}$ MW, minimum capacity $\underline{p}_{i}=\{10,10\}$ MW, adjustment capability $\overline{r}_{i}=\{10,10\}$ MW, and cost coefficients $c_{2i}=\{1,3\}$ and $c_{1i}=\{10,3\}$. The mean forecast of renewable power production is $w^{f}=50$ MW that supplies a half of system load $L=100$ MW. The renewable spillage cost $c^{c}$ and load curtailment cost $c^{s}$ are set to \$100/MWh and \$300/MWh, respectively. We set factor $\rho$ to a small value of $10^{-5}$ in \eqref{algorithm} to isolate forecast impacts from those caused by algorithmic errors. Finally, we notice that convergence time is kept below several minutes. 

\subsection{Impacts of the size of sample-based dataset of producers}\label{exp1}
\begin{figure}
\centering
\resizebox{0.5\textwidth}{!}{%
\begin{tikzpicture}[thick,scale=1]
\begin{axis}[
yshift=4.4cm,
ymin=6,
ymax=240,
xmax=9,
xmin=1,
xtick={1,2,3,4,5,6,7,8,9},
xticklabels={},
typeset ticklabels with strut,
try min ticks=4,
max space between ticks=20pt,
ylabel={$\norm{\fat{\hat{\sigma}}_{1}^{\fat{2}}-\fat{\hat{\sigma}}_{2}^{\fat{2}}}$},
y label style=  {at={(axis description cs:0.1,.5)},rotate=0,anchor=south, font=\scriptsize},
legend style={draw=none, fill=none, legend columns=1, font=\scriptsize, legend pos={north east},xshift=-1.0cm,yshift=-0.05cm},
tick label style={font=\scriptsize},
width=9cm,
height=3.5cm,
]
\addplot [smooth,line width=0.25mm,draw=blue] table [x index = 0,y index = 10] {plots_data/results_ex1_with_comp.dat};\addlegendentry{$\hat{\sigma}_{i}$};
\end{axis}
\begin{axis}[
yshift=4.4cm,
ymin=0,
ymax=62,
xmax=9,
xmin=1,
xtick={1,2,3,4,5,6,7,8,9},
xticklabels={},
typeset ticklabels with strut,
try min ticks=4,
max space between ticks=20pt,
axis y line*=right,
ylabel style = {align=center},
ylabel={$\norm{\fat{\underline{\overline{w}}}_{1}-\fat{\underline{\overline{w}}}_{2}}$},
y label style=  {at={(axis description cs:1.29,.5)},rotate=0,anchor=south, font=\scriptsize},
legend style={draw=none, fill=none, legend columns=1, font=\scriptsize, legend pos={north east},xshift=0.1cm,yshift=-0.04cm},
tick label style={font=\scriptsize},
width=9cm,
height=3.5cm,
]
\addplot [smooth,line width=0.25mm,draw=red] table [x index = 0,y index = 11] {plots_data/results_ex1_with_comp.dat};\addlegendentry{$\underline{\overline{w}}_{i}$};
\end{axis}
\begin{axis}[
yshift=0cm,
ymin=2.25,
ymax=2.8,
xmax=9,
xmin=1,
xtick={1,2,3,4,5,6,7,8,9},
xticklabels={$10$,$30$,$50$,$10^{2}$,$3\hspace{-0.2em}\cdot\hspace{-0.3em}10^2$,$5\hspace{-0.2em}\cdot\hspace{-0.3em}10^2$,$10^{3}$,$1.5\hspace{-0.2em}\cdot\hspace{-0.3em}10^{3}$,$10^{4}$},
typeset ticklabels with strut,
try min ticks=4,
max space between ticks=20pt,
xlabel={The number of samples within the dataset of each producer $(N_{i}^{s})$},
ylabel={$\mathcal{C}$ [$\times10^{3}$]},
x label style = {at={(axis description cs:0.5,-0.18)},anchor=south, font=\scriptsize},
y label style=  {at={(axis description cs:0.1,.5)},rotate=0,anchor=south, font=\scriptsize},
legend style={draw=none, fill=none, legend columns=1, font=\scriptsize, legend pos={north east},xshift=-1.75cm,yshift=-0.04cm},
tick label style={font=\scriptsize},
width=9cm,
height=3.5cm,
]
\addplot [forget plot,smooth,line width=0.01mm,draw=blue!50, name path=mean_min] table [x index = 0,y index = 1] {plots_data/results_ex1_with_comp.dat};
\addplot [forget plot,smooth,line width=0.01mm,draw=blue!50, name path=mean_max] table [x index = 0,y index = 3] {plots_data/results_ex1_with_comp.dat};
\addplot[forget plot,color=blue!50,opacity=0.5]fill between[of=mean_min and mean_max];
\addplot [smooth,draw=blue,line width=0.25mm] table [x index = 0,y index = 2] {plots_data/results_ex1_with_comp.dat};\addlegendentry{$\mathcal{C}$};
\end{axis}
\begin{axis}[
yshift=0cm,
ymin=2,
ymax=15,
xmax=9,
xmin=1,
xtick={1,2,3,4,5,6,7,8,9},
xticklabels={},
typeset ticklabels with strut,
try min ticks=4,
max space between ticks=20pt,
axis y line*=right,
ylabel style = {align=center},
ylabel={$\mathrm{CVaR}_{(5\%)} [\times10^{3}]$},
y label style=  {at={(axis description cs:1.29,.5)},rotate=0,anchor=south, font=\scriptsize},
legend style={draw=none, fill=none, legend columns=1, font=\scriptsize, legend pos={north east},xshift=0.3cm,yshift=-0.04cm},
tick label style={font=\scriptsize},
width=9cm,
height=3.5cm,
]
\addplot [forget plot,smooth,line width=0.01mm,draw=red!50, name path=max_min] table [x index = 0,y index = 4] {plots_data/results_ex1_with_comp.dat};
\addplot [forget plot,smooth,line width=0.01mm,draw=red!50, name path=max_max] table [x index = 0,y index = 6] {plots_data/results_ex1_with_comp.dat};
\addplot[forget plot,color=red!50,opacity=0.5]fill between[of=max_min and max_max];
\addplot [smooth,draw=red,line width=0.25mm] table [x index = 0,y index = 5] {plots_data/results_ex1_with_comp.dat};
\addlegendentry{$\mathrm{CVaR}_{(5\%)}$};
\end{axis}
\begin{axis}[
yshift=2.2cm,
ymin=91,
ymax=100,
xmax=9,
xmin=1,
xtick={1,2,3,4,5,6,7,8,9},
xticklabels={},
ylabel={(1-$\nu$) [\%]},
x label style = {at={(axis description cs:0.5,-0.18)},anchor=south, font=\scriptsize},
y label style=  {at={(axis description cs:0.1,.5)},rotate=0,anchor=south, font=\scriptsize},
legend style={draw=none, fill=none, legend columns=1, font=\scriptsize, legend pos={north west}},
tick label style={font=\scriptsize},
width=9cm,
height=3.5cm,
]
\addplot [smooth,line width=0.01mm,draw=blue!20, name path=mean_min] table [x index = 0,y index = 7] {plots_data/results_ex1_with_comp.dat};
\addplot [smooth,line width=0.01mm,draw=blue!20, name path=mean_max] table [x index = 0,y index = 9] {plots_data/results_ex1_with_comp.dat};
\addplot[color=blue!20]fill between[of=mean_min and mean_max];
\addplot [smooth,draw=blue,line width=0.25mm] table [x index = 0,y index = 8] {plots_data/results_ex1_with_comp.dat};
\end{axis}
\end{tikzpicture}
}
\caption{Impacts of the size of sample-based dataset of producers: The upper plot illustrates the data dissimilarity among producers in the variance ($\hat{\sigma}_{i}$) and interval ($\underline{\overline{w}}_{i}$) of the forecast error distribution. The intermediate plot depicts the empirical reliability level of the system (1-$\nu$). The lower plot shows the average out-of-sample cost ($\mathcal{C}$) and the expectation across 5\% of the worst-case scenarios ($\mathrm{CVaR}_{(5\%)}$). The upper and lower bounds of the envelopes around the solid lines are obtained from the variance of an indicator over 50 simulation runs normalized to the average.}
\label{plots_ex1}
\end{figure}
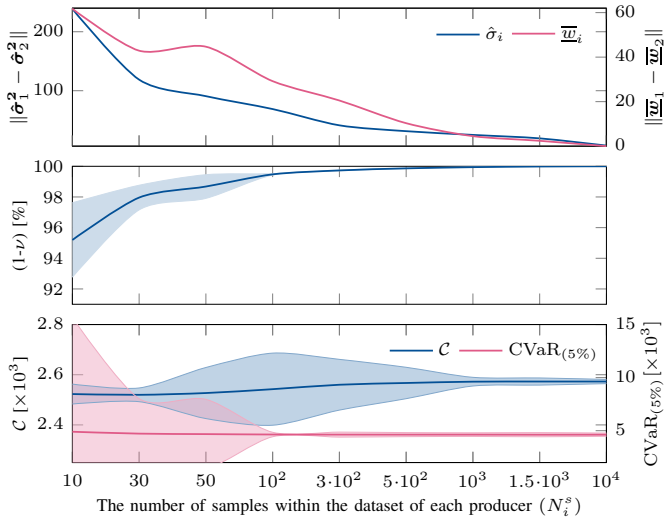

The first experiment assesses the connection between the size of each producer's sample-based dataset, i.e., the number of samples in their dataset, and the overall operational efficiency of the market. Without the loss of generality, we assume that renewable power forecast errors follow normal distribution $\mathcal{N}(0,\sigma^2)$ with $\sigma^{2}=50,$ known to each data provider. The producers, though, receive only $N_{i}^{s}$ number of samples, ranging from $10$ to $10^{4}$. Therefore, producers have different estimates of the variance and forecast error bounds that gradually improve in $N_{i}^{s}$. To measure the level of forecast asymmetry among producers for some $N_{i}^{s}$, we use the $\ell_{2}$-norm as a dissimilarity measure. We consider $\norm{\fat{\hat{\sigma}}^{\fat{2}}_{1}-\fat{\hat{\sigma}}^{\fat{2}}_{2}}$, where $\fat{\hat{\sigma}}^{\fat{2}}_{i}$ is a $1\times50$ vector of variance estimates of producer $i$ obtained after 50 simulations of the evaluation framework. In the same way, we measure the dissimilarity of the interval of forecast error distribution defined as $\fat{\underline{\overline{w}}}_{i}=\fat{\overline{w}}_{i}-\fat{\underline{w}}_{i}.$ 

The projections of the size of sample-based dataset of producers as their private forecast data on the system reliability and out-of-sample cost are depicted in Fig. \ref{plots_ex1}. With data scarcity, the dissimilarity of forecasts among producers is kept high so they have different expectations over underlying uncertainty distribution. As a result, the system runs at a poor expected reliability level with large variance, subsequently resulting in a large variation of the cost in 5\% of the worst-case scenarios ($\mathrm{CVaR}_{(5\%)}$) due to emergent renewable power curtailment and load shedding. This scenario corresponds to the first entry in Fig. \ref{example}.  With increasing size of sample-based dataset, forecasts of producers gradually align, so the variance of the empirical system reliability reduces, and the expected reliability level exceeds 99\% with the size of sample-based dataset $N_{i}^{s}=10^{2}$, yielding the smallest variance of $\mathrm{CVaR}_{(5\%)}$. We relate the range $N_{i}^{s}\in[10^2,10^{3}]$ to the second entry in Fig. \ref{example}. 
Notice, the equilibrium solution obtained with  $N_{i}^{s}\geqslant10^3$ converges to the ideal solution provided by the centralized model \eqref{CC_dispatch} using perfect forecast. Finally, we notice that the empirical reliability level in Fig. \ref{plots_ex1} serves as a proxy function to the condition on size of sample-based dataset in \cite{margellos2014road} that can be used to imposed probabilistic guarantee for the decentralized algorithm in \eqref{algorithm}. 

\subsection{Learning from data}\label{exp2}
\begin{figure}
\centering
\resizebox{0.47\textwidth}{!}{%
\begin{tikzpicture}[thick,scale=1]
\begin{axis}[
yshift=4.4cm,
ymin=0,
ymax=100,
xmax=10,
xmin=1,
xtick={1,2,3,4,5,6,7,8,9,10},
xticklabels={},
typeset ticklabels with strut,
try min ticks=4,
max space between ticks=20pt,
ylabel={$\norm{\circ_{1}-\circ_{2}}$},
y label style=  {at={(axis description cs:0.1,.5)},rotate=0,anchor=south, font=\scriptsize},
legend style={draw=none, fill=none, legend columns=2, font=\scriptsize, legend pos={north east},yshift=-0.0cm},
tick label style={font=\scriptsize},
width=9cm,
height=3.5cm,
]
\addplot [smooth,line width=0.25mm,draw=red] table [x index = 0,y index = 19] {plots_data/results_ex2_final.dat};
\addlegendentry{${\hat{\alpha}}$};
\addplot [smooth,line width=0.25mm,draw=blue] table [x index = 0,y index = 20] {plots_data/results_ex2_final.dat};
\addlegendentry{${\hat{\beta}}$};
\end{axis}
\begin{axis}[
yshift=0cm,
ymin=2.25,
ymax=2.95,
xmax=10,
xmin=1,
xtick={1,2,3,4,5,6,7,8,9,10},
xticklabels={$10$,$20$,$30$,$40$,$50$,$10^{2}$,$3\hspace{-0.2em}\cdot\hspace{-0.3em}10^2$,$5\hspace{-0.2em}\cdot\hspace{-0.3em}10^2$,$1\hspace{-0.3em}\cdot\hspace{-0.3em}10^{3}$,$1.5\hspace{-0.2em}\cdot\hspace{-0.3em}10^{3}$},
typeset ticklabels with strut,
try min ticks=4,
max space between ticks=20pt,
xlabel={The number of samples within the dataset of each producer $(N_{i}^{s})$},
ylabel={$\mathcal{C}$ [$\times10^{3}$]},
x label style = {at={(axis description cs:0.5,-0.18)},anchor=south, font=\scriptsize},
y label style=  {at={(axis description cs:0.1,.5)},rotate=0,anchor=south, font=\scriptsize},
legend style={draw=none, fill=none, legend columns=2, font=\scriptsize, legend pos={north east},yshift=-0.0cm},
tick label style={font=\scriptsize},
width=9cm,
height=3.5cm,
]
\addplot [forget plot,smooth,line width=0.01mm,draw=red!50, name path=mean_min] table [x index = 0,y index = 10] {plots_data/results_ex2_final.dat};
\addplot [forget plot,smooth,line width=0.01mm,draw=red!50, name path=mean_max] table [x index = 0,y index = 12] {plots_data/results_ex2_final.dat};
\addplot[forget plot,color=red!50,opacity=0.5]fill between[of=mean_min and mean_max];
\addplot [smooth,draw=red,line width=0.25mm] table [x index = 0,y index = 11] {plots_data/results_ex2_final.dat};\addlegendentry{no learning};
\addplot [forget plot,smooth,line width=0.01mm,draw=blue!50, name path=mean_min] table [x index = 0,y index = 1] {plots_data/results_ex2_final.dat};
\addplot [forget plot,smooth,line width=0.01mm,draw=blue!50, name path=mean_max] table [x index = 0,y index = 3] {plots_data/results_ex2_final.dat};
\addplot[forget plot,color=blue!50,opacity=0.5]fill between[of=mean_min and mean_max];
\addplot [smooth,draw=blue,line width=0.25mm] table [x index = 0,y index = 2] {plots_data/results_ex2_final.dat};\addlegendentry{learning};
\end{axis}
\begin{axis}[
yshift=2.2cm,
ymin=91,
ymax=100,
xmax=10,
xmin=1,
xtick={1,2,3,4,5,6,7,8,9,10},
xticklabels={},
ylabel={(1-$\nu$) [\%]},
x label style = {at={(axis description cs:0.5,-0.18)},anchor=south, font=\scriptsize},
y label style=  {at={(axis description cs:0.1,.5)},rotate=0,anchor=south, font=\scriptsize},
legend style={draw=none, fill=none, legend columns=2, font=\scriptsize, legend pos={north east},yshift=-5},
tick label style={font=\scriptsize},
width=9cm,
height=3.5cm,
]
\addplot [forget plot,smooth,line width=0.01mm,draw=red!20, name path=mean_min] table [x index = 0,y index = 16] {plots_data/results_ex2_final.dat};
\addplot [forget plot,smooth,line width=0.01mm,draw=red!20, name path=mean_max] table [x index = 0,y index = 18] {plots_data/results_ex2_final.dat};
\addplot[forget plot,color=red!50,opacity=0.5]fill between[of=mean_min and mean_max];
\addplot [smooth,draw=red,line width=0.25mm] table [x index = 0,y index = 17] {plots_data/results_ex2_final.dat};
\addlegendentry{no learning};
\addplot [forget plot,smooth,line width=0.01mm,draw=blue!20, name path=mean_min] table [x index = 0,y index = 7] {plots_data/results_ex2_final.dat};
\addplot [forget plot,smooth,line width=0.01mm,draw=blue!20, name path=mean_max] table [x index = 0,y index = 9] {plots_data/results_ex2_final.dat};
\addplot[forget plot,color=blue!50,opacity=0.5]fill between[of=mean_min and mean_max];
\addplot [smooth,draw=blue,line width=0.25mm] table [x index = 0,y index = 8] {plots_data/results_ex2_final.dat};
\addlegendentry{learning};
\end{axis}
\end{tikzpicture}
}
\caption{Learning from data: The upper plot shows the data dissimilarity of producers in the estimators of beta distribution $(\hat{\alpha}_{i},\hat{\beta}_{i})$. The intermediate plot illustrates the empirical system reliability (1-$\nu$). The lower plot depicts the average out-of-sample cost ($\mathcal{C}$). The upper and lower bounds of the envelopes around the solid lines are obtained from the variance of an indicator over 50 simulation runs normalized to the average.}
\label{plots_ex2}
\end{figure}
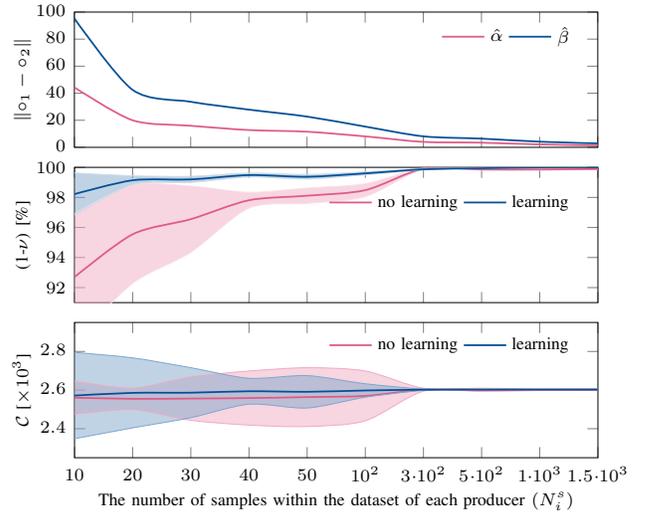

When data is limited, it is reasonable to assume that producers learn statistical properties of underlying uncertainty from the data at hand. In this experiment, we consider that the renewable power forecast errors are obtained from a beta distribution with parameters $\tilde{\alpha}=5$ and $\tilde{\beta}=10$ scaled by a factor of 65. We assume that the producers are aware of the type of distribution and estimate their parameters from datasets $\mathcal{D}_{i}$ relying on maximum likelihood estimation. Once producer estimates $\hat{\alpha}_{i}$ and $\hat{\beta}_{i}$ are obtained, they generate $N=3\cdot10^5$ number of samples in attempt to recover the true distribution. Then, each producer enriches the initial dataset $\mathcal{D}_{i}$ with generated $N$ samples and enforces constraints \eqref{reformulated_cc} over a new set of samples. Similarly to the previous experiment, we consider $\ell_{2}$-norm as a dissimilarity measure for $\hat{\alpha}_{i}$ and $\hat{\beta}_{i}$.

The results of this experiment are summarized in Fig. \ref{plots_ex2}. With a small number of samples, the two producers demonstrate highly divergent estimates of the true parameters of beta distribution. However, the system reliability significantly improves compared to the direct implementation of the data at hand. This observation is aligned with the randomization approach in \cite{margellos2014road}. On the other hand, the poor statistical estimates of the true distribution yield the large variance of the expected operating cost. The increasing size of sample-based dataset improves the statistical significance of producers estimates. Eventually, when producers data is subject to learning,  the system reaches a high level of reliability and small variance of operating cost significantly sooner than in the reference case with no learning.

\subsection{Data sharing}\label{exp3}
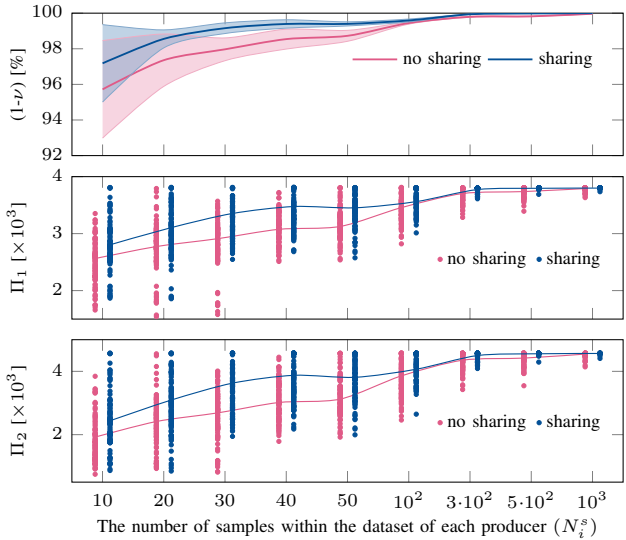
\begin{figure}
\centering
\resizebox{0.47\textwidth}{!}{%
\begin{tikzpicture}[thick,scale=1]
\begin{axis}[
every axis plot/.append style=only marks,
yshift=-2.2cm,
ymin=1.5,
ymax=4.0,
xmax=9.5,
xmin=0.5,
xtick={1,2,3,4,5,6,7,8,9},
xticklabels={},
typeset ticklabels with strut,
try min ticks=4,
max space between ticks=20pt,
ylabel={$\Pi_{1}$ [$\times10^{3}$]},
y label style={at={(axis description cs:0.1,.5)},rotate=0,anchor=south, font=\scriptsize},
x label style={at={(axis description cs:0.5,-0.18)},anchor=south, font=\scriptsize},
tick label style={font=\scriptsize},
width=9cm,
height=3.5cm,
]
\end{axis}
\begin{axis}[
yshift=-2.2cm,
every axis plot/.append style=only marks,
axis line style={draw=none},
tick style={draw=none},
xshift=-0.1cm,
ymin=1.5,
ymax=4.0,
xmax=9.5,
xmin=0.5,
xtick={1,2,3,4,5,6,7,8,9},
xticklabels={},
yticklabels={},
ylabel={},
xlabel={},
legend style={draw=none, fill=none, legend columns=2, font=\scriptsize, legend pos={north east},yshift=-0.8cm,xshift=-0.8cm},
width=9cm,
height=3.5cm,
]
\addplot+[color=red,only marks,mark=*,mark options={solid},mark size = 0.08em] table [x index = 0,y index = 1] {plots_data/results_ex3_profits_v2.dat};\addlegendentry{no sharing};
\addplot [smooth,draw=red,line width=0.15mm] table [x index = 0,y index = 8] {plots_data/results_ex3.dat};
\end{axis}
\begin{axis}[
yshift=-2.2cm,
every axis plot/.append style=only marks,
axis line style={draw=none},
tick style={draw=none},
xshift=0.1cm,
ymin=1.5,
ymax=4.0,
xmax=9.5,
xmin=0.5,
xtick={1,2,3,4,5,6,7,8,9},
xticklabels={},
yticklabels={},
ylabel={},
xlabel={},
legend style={draw=none, fill=none, legend columns=2, font=\scriptsize, legend pos={north east},yshift=-0.8cm},
width=9cm,
height=3.5cm,
]
\addplot+[color=blue,only marks,mark=*,mark options={solid},mark size = 0.08em] table [x index = 0,y index = 3] {plots_data/results_ex3_profits_v2.dat};\addlegendentry{sharing};
\addplot [smooth,draw=blue,line width=0.15mm] table [x index = 0,y index = 11] {plots_data/results_ex3.dat};
\end{axis}
\begin{axis}[
every axis plot/.append style=only marks,
yshift=-4.4cm,
ymin=0.5,
ymax=5,
xmax=9.5,
xmin=0.5,
xtick={1,2,3,4,5,6,7,8,9},
xticklabels={$10$,$20$,$30$,$40$,$50$,$10^{2}$,$3\hspace{-0.2em}\cdot\hspace{-0.3em}10^2$,$5\hspace{-0.2em}\cdot\hspace{-0.3em}10^2$,$10^3$},
typeset ticklabels with strut,
try min ticks=4,
max space between ticks=20pt,
ylabel={$\Pi_{2}$ [$\times10^{3}$]},
xlabel={The number of samples within the dataset of each producer $(N_{i}^{s})$},
y label style={at={(axis description cs:0.1,.5)},rotate=0,anchor=south, font=\scriptsize},
x label style={at={(axis description cs:0.5,-0.18)},anchor=south, font=\scriptsize},
tick label style={font=\scriptsize},
width=9cm,
height=3.5cm,
]
\end{axis}
\begin{axis}[
yshift=-4.4cm,
every axis plot/.append style=only marks,
axis line style={draw=none},
tick style={draw=none},
xshift=-0.1cm,
ymin=0.5,
ymax=5,
xmax=9.5,
xmin=0.5,
xtick={1,2,3,4,5,6,7,8,9},
xticklabels={},
yticklabels={},
ylabel={},
xlabel={},
legend style={draw=none, fill=none, legend columns=2, font=\scriptsize, legend pos={north east},yshift=-0.8cm,xshift=-0.8cm},
width=9cm,
height=3.5cm,
]
\addplot+[color=red,only marks,mark=*,mark options={solid},mark size = 0.08em] table [x index = 0,y index = 2] {plots_data/results_ex3_profits_v2.dat};\addlegendentry{no sharing};
\addplot [smooth,draw=red,line width=0.15mm] table [x index = 0,y index = 14] {plots_data/results_ex3.dat};
\end{axis}
\begin{axis}[
yshift=-4.4cm,
every axis plot/.append style=only marks,
axis line style={draw=none},
tick style={draw=none},
xshift=0.1cm,
ymin=0.5,
ymax=5,
xmax=9.5,
xmin=0.5,
xtick={1,2,3,4,5,6,7,8,9},
xticklabels={},
yticklabels={},
ylabel={},
xlabel={},
legend style={draw=none, fill=none, legend columns=2, font=\scriptsize, legend pos={north east},yshift=-0.8cm},
width=9cm,
height=3.5cm,
]
\addplot+[color=blue,only marks,mark=*,mark options={solid},mark size = 0.08em] table [x index = 0,y index = 4] {plots_data/results_ex3_profits_v2.dat};\addlegendentry{sharing};
\addplot [smooth,draw=blue,line width=0.15mm] table [x index = 0,y index = 17] {plots_data/results_ex3.dat};
\end{axis}
\begin{axis}[
yshift=0cm,
ymin=92,
ymax=100,
xmax=9.5,
xmin=0.5,
xtick={1,2,3,4,5,6,7,8,9},
xticklabels={},
typeset ticklabels with strut,
try min ticks=4,
max space between ticks=20pt,
ylabel={(1-$\nu$) [\%]},
x label style = {at={(axis description cs:0.5,-0.18)},anchor=south, font=\scriptsize},
y label style=  {at={(axis description cs:0.1,.5)},rotate=0,anchor=south, font=\scriptsize},
legend style={draw=none, fill=none, legend columns=2, font=\scriptsize, legend pos={north east},yshift=-0.3cm},
tick label style={font=\scriptsize},
width=9cm,
height=3.5cm,
]

\addplot [forget plot,smooth,line width=0.01mm,draw=red!50, name path=mean_min] table [x index = 0,y index = 1] {plots_data/results_ex3.dat};
\addplot [forget plot,smooth,line width=0.01mm,draw=red!50, name path=mean_max] table [x index = 0,y index = 3] {plots_data/results_ex3.dat};
\addplot[forget plot,color=red!50,opacity=0.5]fill between[of=mean_min and mean_max];
\addplot [smooth,draw=red,line width=0.25mm] table [x index = 0,y index = 2] {plots_data/results_ex3.dat};\addlegendentry{no sharing};

\addplot [forget plot,smooth,line width=0.01mm,draw=blue!50, name path=mean_min] table [x index = 0,y index = 4] {plots_data/results_ex3.dat};
\addplot [forget plot,smooth,line width=0.01mm,draw=blue!50, name path=mean_max] table [x index = 0,y index = 6] {plots_data/results_ex3.dat};
\addplot[forget plot,color=blue!50,opacity=0.5]fill between[of=mean_min and mean_max];
\addplot [smooth,draw=blue,line width=0.25mm] table [x index = 0,y index = 5] {plots_data/results_ex3.dat};\addlegendentry{sharing};
\end{axis}
\end{tikzpicture}
}
\caption{Data sharing: The upper plot shows the average and variance of empirical system reliability (1-$\nu$) over 50 simulation runs. The intermediate and lower plots depict the expected payoffs ($\Pi_{1}$ and $\Pi_{2}$) of the two producers in 50 simulation runs, where the solid line illustrates the average value.}
\label{plots_ex3}
\end{figure}

Data sharing among producers can be seen as a natural way to enhance operational performance of the system. However, the data sharing among producers is meaningful if producers have strong incentives. To estimate producer incentives, we consider the payoff of each generator in some admissible uncertainty scenario $s$ according to the following function
\begingroup
\allowdisplaybreaks
\begin{align*}
    \Pi_{is} = \lambda^{\text{e}}p_{i} + &\lambda^{\text{r}}\alpha_{i} - c_{2i}\braceit{p_{i}+ r_{is}}^2 - c_{1i}\braceit{p_{i} + r_{is}},
\end{align*}
\endgroup
where the first two terms compute the revenue by selling energy and reserve, while the last two terms define the actual cost incurred under scenario $s$. 

Similarly to the first experiment, we consider forecast error distribution $\mathcal{N}(0,\sigma^2)$ with $\sigma^{2}=50$. The system reliability is opposed to the expected payoff of producers for different size of sample-based datasets in Fig. \ref{plots_ex3}. We observe that up to size of sample-based dataset $N_{i}^{s}=3\cdot10^{2}$, data sharing significantly improves the system reliability, but loses its value with a larger number of samples as the dataset of producers becomes sufficiently large. The improvement in system reliability, however, does not come at the expense of any of  producers. In contrast, we observe that the average expected payoff for the two producers is higher and with smaller variation when they  share data. Moreover, the lower bound on  payoffs tends to increase with data sharing for almost any  size of sample-based dataset. 

\section{Conclusion}

This paper models the renewable power forecast asymmetry among market participants and studies its impacts on market equilibrium. The forecast of each participant is modeled in form of private samples, internalized into their profit-maximizing objectives and constraint sets. Our numerical experiments show that the market-clearing outcomes are strongly conditioned by the quality and asymmetry of private forecasts. The system reliability and operating cost improve significantly with the increasing availability of historical renewable power observations. However, if forecast availability is limited, the system and the producers may individually benefit from learning statistical properties of the data at hand and from sharing their private forecasts.

\balance
\bibliographystyle{IEEEtran}
\bibliography{refs}
\end{document}